\documentclass[12pt,final]{elsarticle}

\usepackage{amssymb}
\usepackage{graphicx, color}
\usepackage{esvect}
\usepackage{rawfonts}
\usepackage{pictexwd}
\usepackage{subcaption}
\usepackage{physics}
\usepackage{mhchem}
\usepackage[utf8]{inputenc}
\usepackage[T1]{fontenc}
\usepackage{cancel}
\usepackage{algorithm}
\usepackage{algpseudocode}
\usepackage{booktabs}
\usepackage{multirow}

\journal{}

\begin{document}

\begin{frontmatter}

\title{Application of a Temporal Multiscale Method for Efficient Simulation of Degradation in PEM Water Electrolysis under Dynamic Operation}

 \author[label1]{Dayron Chang Dominguez}
 \author[label2]{An Phuc Dam}
 \author[label1]{Thomas Richter}
 \author[label2]{Kai Sundmacher}
 \author[label3]{Shaun M. Alia}

 \affiliation[label1]{organization={Otto-von-Guericke University},
             addressline={Universitätsplatz 2},
             city={Magdeburg},
             postcode={39106},
             state={Saxony-Anhalt},
             country={Germany}}

 \affiliation[label2]{organization={Max-Planck-Institute for Dynamics of Complex Technical Systems},
             addressline={SandtorStraße 1},
             city={Magdeburg},
             postcode={39106},
             state={Saxony-Anhalt},
             country={Germany}}
             
\affiliation[label3]{
			 organization={National Renewable Energy Laboratory},
             addressline={15013 Denver West Parkway},
             city={Golden},
             postcode={CO 80401},
             state={Colorado},
             country={United States of America}
			}

\begin{abstract}

Hydrogen is vital for sectors like chemicals and others, driven by the need to reduce carbon emissions. Proton Electrolyte Membrane Water Electrolysis (PEMWE) is a key technology for the production of green hydrogen under fluctuating conditions of renewable power sources. However, due to the scarcity of noble metal materials, the stability of the anode catalyst layer under dynamic operating conditions must be better understood. Model-aided investigation approaches are essential due to the back-box nature of the electrochemical system and the high costs of experimental long-term testing. In this work, a temporal multiscale method based on a Heterogeneous technique is applied to reduce the computational effort of simulating long-term degradation, focused on catalyst dissolution. Such an approach characterizes the problem in fast locally periodic processes, influenced by the dynamic operation and slow processes attributed to the gradual degradation of the catalyst layer. A mechanistic model that includes the oxygen evolution reaction, catalyst dissolution and hydrogen permeation from the cathode to the anode side is hypothesized and implemented. The multiscale approach notably reduces computational effort of simulation from hours to mere minutes. This efficiency gain is ascribed to the limited evolution of Slow-Scale variables during each period of time of the Fast-Scale variables. Consequently, simulation of the fast processes is required only until local periodicity is achieved within each Slow-Scale time step. Thus, the developed temporal multiscale approach proves to be highly effective in accelerating parameter estimation and predictive simulation steps, as could be verified through the results of this article. In this way, the method can support systematic model development to describe degradation in PEMWE under dynamic operating conditions.
\end{abstract}


\end{frontmatter}

\section*{Introduction}
Degradation of the anode catalyst layer is a major challenge in Proton Exchange Membrane Water Electrolysis (PEMWE), especially due to the need to reduce loadings of noble metal materials such as iridium. Understanding the degradation phenomena that occur within the ACL of the PEM Water Electrolyzer \cite{prestat_corrosion_2023,Dam2023} is challenging due to the “black-box” nature of such a PEMWE cell during operation, even more so under dynamic load conditions \cite{voronova_effect_2022}. Insights are typically gained by studying the system in- and output behavior as well as by conducting post-test analysis of the system. 

In situ experimental knowledge regarding the phenomenon of Ir catalyst dissolution has been acquired through the use of so-called half-cell setups, where liquid electrolytes are employed, and models are validated \cite{Qian2022}. However, obtaining similar in situ data is challenging for the industrially relevant membrane electrode assembly (MEA) setups. Transfer of knowledge of catalyst dissolution from half-cells to full-cells is not straight forward \cite{dam_role_2020, knoppel_limitations_2021}.

Therefore, mathematical modelling is of essential importance to undestand the physico-chemical processes occuring inside the cell. Since the early 2000s modelling processes in PEMWEs has been used as a tool to make up for shortfalls in such a closed schema \cite{grigorev_electrolysis_2001, onda_performance_2002}. An extensive classification of different ways of modelling PEMWE can be found in \cite{olivier_low-temperature_2017}. Some are focused on the understanding of the phenomena under various conditions such as high-pressure water electrolysis \cite{grigoriev_mathematical_2010, yigit_mathematical_2016} or hydrogen permeation \cite{trinke_hydrogen_2016, papakonstantinou_h2_2019, franz_transient_2023}. Additionally, models ranging from 0D to 3D of PEMWEs have been developed based on physical laws to comprehend processes such as mass transport \cite{abdol_rahim_overview_2016,aubras_two-dimensional_2017}, fluid behavior through the Porous Transport Layer (PTL) \cite{ojong_development_2017,doan_review_2021,lin_1d_2022,ni_anode_2023} among others.

In the area of performance characterization and control system development, there are works like \cite{gorgun_dynamic_2006,Marangio2009,lebbal_identification_2009, espinosa-lopez_modelling_2018,ogumerem_parametric_2020,Dizon2021}, mainly focus on analysis of experimental results and long-term durability evaluation. The authors utilized fitting techniques based on experimental data to identify parameters and optimize the performance of PEMWEs using often “black-box” models. In this regard, the machine learning era unfolds data-driven simulating techniques which accelerate the fitting processes due to its computational efficiency \cite{ding_guiding_2022}.

\begin{figure}[h]
    \begin{center}
	  \includegraphics[width=1.0\textwidth]{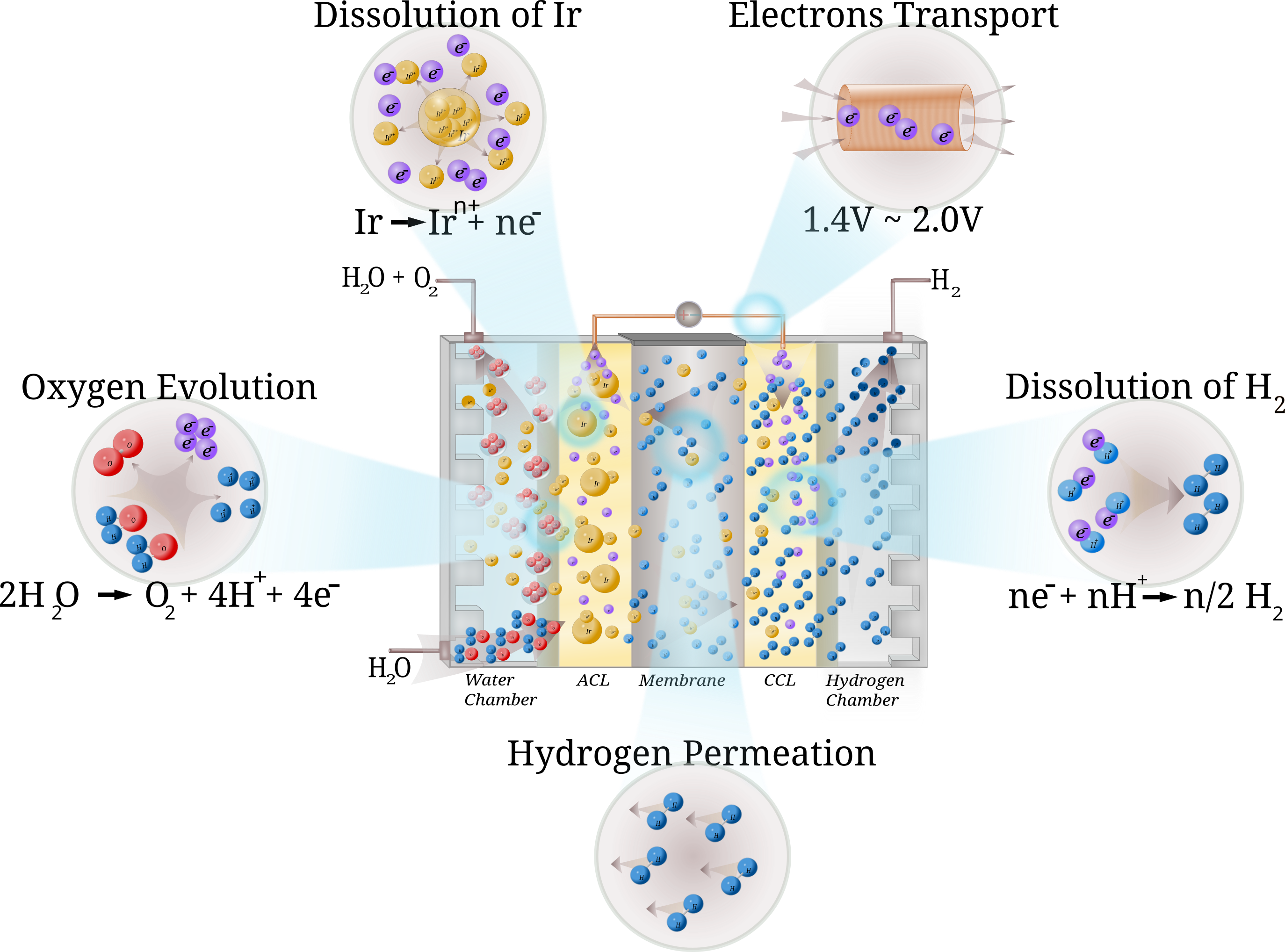}
    \end{center}
    \caption{Theoretical schema of the processes of the PEMWE that are taken into consideration in the proposed model.}
    \label{fig:PEMWE_schema}
\end{figure}

Models have been developed for PEMWE in which mutliscale analysis has applied in order to better resolve spatial dimensions  \cite{Franco2007, Franco2008, Franco2009, Oliveira2012}. In this study a white-box methodology for modelling the dissolution phenomena inside the PEMWE ACL shall be developed (Figure \ref{fig:PEMWE_schema}). This will be numerically supported by applying a temporal multiscale method for the time dimension in order to reduce computational effort.

This technique has been used to describe long-term effects in the interaction of processes of different time scales in different applications as is shown in \cite{Yang2016,Brinkmann2018}. In this way, the Heterogeneous Multiscale Method (HMM) \cite{Weinan2003} is one of the most promising options to efficiently decouple the micro-scale and macro-scale of the problem solving the last one based on temporal averaging of the first one. A review of these techniques can be found in \cite{e_heterogeneous_2007, gravemeier_towards_2008} where authors compare the HMM  with other techniques, mainly applied for ODEs or to spatial multiscales. 

In \cite{frei_efficient_2020} a description of a computational temporal multiscale method is derived. It assumes that the underlying fast process has a time quasi-locally periodic solution if the influence of the Slow-Scale solution is fixed during the slow time step. Therefore, only solving the Fast-Scale variables until such a time quasi-locally periodic is achieved allows the Slow-Scale variables to have very large time steps and accelerate the simulations in about $10^4$ times.

In this paper a simple physico-chemical model of Iridium catalyst dissolution is proposed in Section \ref{sec:1}. It is shown that with the introduction of a numerical temporal multiscale method in Section \ref{sec:2}, the simulation time of the PEMWE ACL long-duration stability tests is drastically reduced. In Section \ref{sec:3} the method is numerically analyzed showing how to overcome the challenges of the stiffness nature of the problem. The model parameters are fitted to experimental data of two dynamic operation profiles \cite{alia_electrolyzer_2019} of accelerated stress tests (see section \ref{sec:4}). A good fit of model simulation with experimental data is achieved. The model with the fitted parameters is then used to predict other dynamic operation profiles. The simulation is compared to experimental data and implications on the hypothesized electrochemical model are discussed.  

Overall, the aim of this work is to demonstrate how the introduction of a temporal multiscale approach, facilitates the systematic development of a mathematical model capable of describing the effects of different dynamic operational profiles on catalyst degradation in PEMWE.

\section{Anode Catalyst Layer degradation model}
\label{sec:1}
In this section, the electrochemical degradation model of the ACL during the operation of the PEMWE is described. The dissolution model is presented in Section \ref{sec:1.1}, and an important influence of hydrogen is assumed, as it can chemically reduce the Ir oxide, thus leading to its destabilization during dynamic operation of PEMWE \cite{papakonstantinou_h2_2019, weis_impact_2019}. Therefore, a transport model has been developed and implemented to describe the dynamic permeation process of hydrogen from the cathode to the anode side which will be discussed in section \ref{sec:1.2}.
\subsection{Iridium dissolution model}\label{sec:1.1}
To describe the dissolution process of Ir, based on knowledge and understanding from previous works in literature, a simple mechanistic model is hypothesized and implemented.
\begin{figure}[h]
    \begin{center}
	  \includegraphics[width=0.8\textwidth]{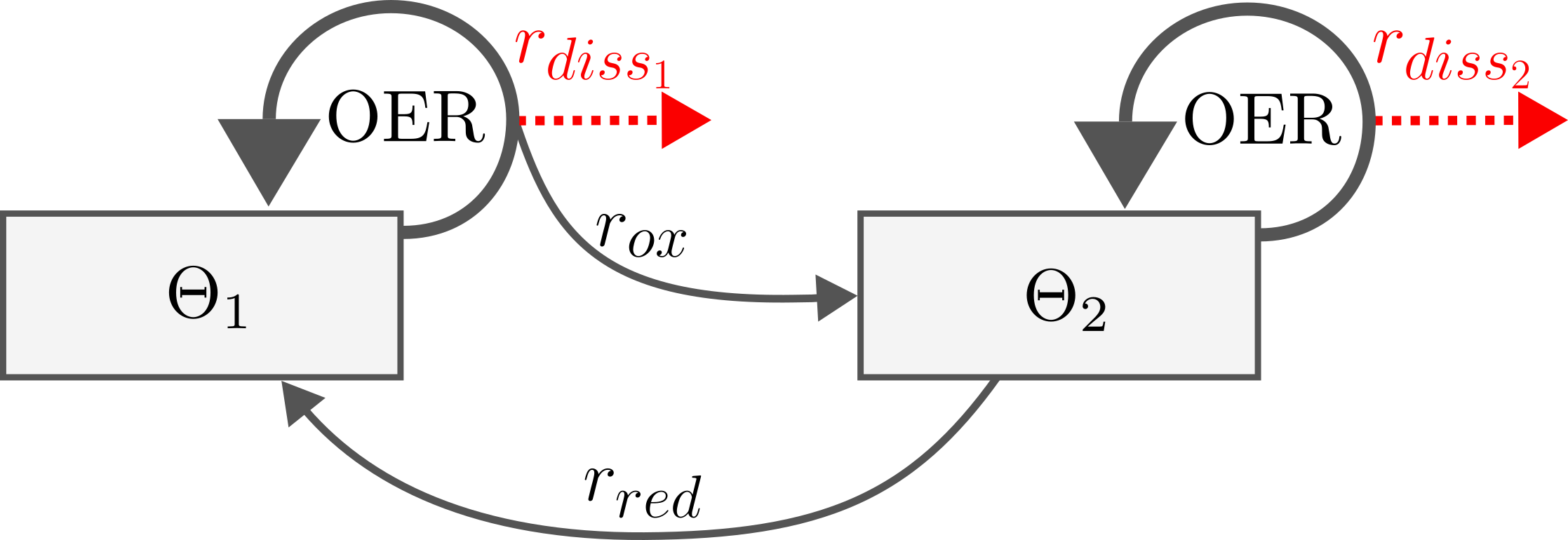}
    \end{center}
    \caption{Theoretical representation of the electrochemical model for corrosion of the material in the ACL. }
    \label{fig:Electrochemical_model}
\end{figure}
The model considers nano-particles of Iridium in which the catalyst material surface can be in two different chemical oxidation states as shown in Figure \ref{fig:Electrochemical_model}. They can be either in a reduced or an oxidized state. $\Theta_1$ and $\Theta_2$ correspond to the more reduced and the more oxidized state, respectively.
	
Figure \ref{fig:Electrochemical_model} illustrates the schematic of Iridium dissolution model which was implemented in this work. It is well known that catalyst dissolution and OER are closely coupled phenomena. An unstable state of the catalyst surface during the OER reaction cycle leads to catalyst dissolution \cite{kasian_common_2018, Geiger2018, dam_role_2020}. 

In this model it is assumed that the surface is generally in an oxidized state during the considered operation of PEMWE at OER potentials. However, hydrogen which permeates from the cathode to the anode side can reduce the oxidation state of the active site on the catalyst surface. Thus, this reduction converts it from the state $\Theta_2$ to the state $\Theta_1$ (`$r_ {red}$'). 

Nevertheless, the OER leads to dissolution of those chemically reduced active sites. As a result, the Ir oxide lattice stucture underneath the dissolved catalyst species becomes exposed to the surface. The surface state $\Theta_1$ then reverts back towards the more stable state $\Theta_2$ (`$r_{ox}$'). However, the dissolution of catalyst species leads to a loss of catalyst material leading to a loss of Electrochemical Surface Area (ECSA) and a gradual decline in electrolyzer performance over time. 

The catalyst surface in state $\Theta_2$ is more stable than in state $\Theta_1$. However, due to the inherent coupling of catalyst dissolution and OER \cite{binninger_thermodynamic_2015, kasian_common_2018}, the catalyst can dissolved also in the former state. Therefore, two dissolution pathways are considered in the implemented model, as shown in Figure \ref{fig:Electrochemical_model}.

The involvement of hydrogen in the dissolution of catalyst material may influence the observed ACL degradation rates.This is especially notable under different dynamic operation conditions, as these conditions strongly influence the dynamic behavior of oxygen and hydrogen partial pressure in the ACL. A high current density results in a low hydrogen partial pressure as evolved oxygen generates a flushing effect on the ACL. On the other hand, the decrease in current density leads to a dynamic accumulation of hydrogen in the ACL, attributable to the prior storage of permeating hydrogen within the membrane.

The implemented mathematical model can be expressed by the following system of ordinary differential (ODE) equations: 
\begin{align}
\dfrac{d \Theta_1}{dt} &= 2 k_r (1- \Theta_1)^2 \left( C^{\text{ACL}}_{H_2} R T \right) - \frac{\Theta_1}{\gamma} k_{\text{diss}_1} e^{f E^{\text{ACL}}}\label{eq:kineticks}\\
\dfrac{d C^{\text{ACL}}_{O_2}}{dt} &= V_{\text{ACL}}^{-1}\left[ \frac{i}{z_{O_2} F} - A_{geo} F^{PTL}_{O_2} \right]\label{eq:balance_O}\\
\dfrac{d C^{\text{ACL}}_{H_2}}{dt} &= \frac{A_{geo}}{V_{\text{ACL}}}\left[ F^{\text{Mem}}_{H_2} - F^{\text{PTL}}_{H_2}\right]\label{eq:balance_H}\\
\dfrac{d N^{np}_{Ir}}{dt} &= - \frac{A^{\text{act}}_{\ce{Ir}}(N^{np}_{Ir}(t))}{A_{geo}} \left[\Theta_1 k_{\text{diss}_1} + \left( 1 - \Theta_1 \right) k_{\text{diss}_2} \right] e^{f E^{\text{ACL}}}\label{eq:dissolution}
\end{align}
The origin and meaning of the variables and terms in this system of ODEs shall be described in the following sections. $\Theta_1$ describes the fraction of active sites which are in state number 1 (see \ref{fig:Electrochemical_model}). The first term of \eqref{eq:kineticks} expresses the reduction rate with its reaction rate constant $k_r$. 

The reduction kinetics depends on the pressure of hydrogen in the ACL since it occurs via chemical reduction with hydrogen and is written via elementary reaction kinetic formulation. The modelled reduction process may be described by the chemical reaction equation $2 IrO_2 \ + \ H_2 \ \rightarrow  2 HIrO_2$. The hydrogen partial pressure can be calculated using the ideal gas law:
\begin{equation}\label{eq:idealGasLaw}
P_{H_2}^{\text{ACL}} = C^{\text{ACL}}_{H_2} R T,
\end{equation}
with $R$ and $T$ being the universal gas constant and the temperature, respectively. $C^{\text{ACL}}_{H_2}$ is the concentration of hydrogen in the ACL and is given in $mol/m^3$. The hydrogen-induced reduction reaction leads to a conversion of active sites from state $\Theta_2$ to the state $\Theta_1$. 

On the other hand, the right term of equation \eqref{eq:kineticks} represents the conversion of sites in state $\Theta_1$ back to the state $\Theta_2$. This is explained by the dissolution of reduced and destabilized sites, exposing the Ir oxide species that were originally located underneath the dissolving sites. $k_{\text{diss}_1}$ is the dissolution rate constant for the sites in $\Theta_1$ state. $\gamma$ is the number of active sites per surface area of the Ir oxide nano-particles and a value of $1.121 \cdot 10^{-4} \, mol/m^2$ is used for the simulation \cite{Papakonstantinou2022}. $E^{\text{ACL}}$ is the potential applied on the ACL, and $f = F/(R T)$, where $F$ is the Faraday constant $ 96485 \ C / mol$.

In terms of spatial direction, the ACL is modelled zero-dimensionally. The balance of the gases it is shown in the equations \eqref{eq:balance_O} and \eqref{eq:balance_H}. The $\ce{O2}$ concentration is balanced by a source and a flow term. Within the source term, $i$ represents the current density ($A/m^2$) and $z_{\ce{O2}}$ is equal to $4$, corresponding with the number of electrons transferred per evolved oxygen molecule. The other term of equation \eqref{eq:balance_O} defines the flow of $\ce{O2}$ from the ACL to the porous space inside the adjacent PTL domain. These fluxes can be described by:
\begin{equation}
F^{\text{PTL}}_{\ce{O2}} = F^{\text{tot}}  X_{\ce{O2}} \quad \text{and} \quad F^{\text{PTL}}_{\ce{H2}} = F^{\text{tot}}  X_{\ce{H2}},\label{eq:ptl_flow_o2_h2}
\end{equation}
where $F^{\text{tot}}$ is the total gas flow out of the ACL and $X_{\ce{O2}}$ and $X_{\ce{H2}}$ are the oxygen and the hydrogen fraction in the gas phase, respectively. These fractions can be calculated by:
\begin{equation}
X_{\ce{O2}} = \frac{P^{\text{ACL}}_{\ce{O2}}}{P^{\text{ACL}}_{\ce{H2}} + P^{\text{ACL}}_{\ce{O2}}} \quad \text{and} \quad X_{\ce{H2}} = \frac{P^{\text{ACL}}_{\ce{H2}}}{P^{\text{ACL}}_{\ce{H2}} + P^{\text{ACL}}_{\ce{O2}}}.\label{eq:fractions_pressure}
\end{equation}
The total gas flow $F^{\text{tot}}$ is calculated using the following equation:
\begin{equation}
F^{\text{tot}} = k_{\text{ACL}} \left(\left( P^{\text{ACL}}_{\ce{O2}} + P^{\text{ACL}}_{\ce{H2}} \right) - P^{\text{atm}}\right)\label{eq:total_flow}.
\end{equation}
The pressure difference serves as the mechanical driving force for the gas flow out of the ACL. The parameter $k_{\text{ACL}}$ is a mass transport coefficient and its inverse describes the resistance of gas transport out of the ACL domain. Although, in reality, this transport resistance is distributed spatially across the ACL thickness, for the sake of simplicity, this work assumes discrete occurrence at the interface between the ACL and PTL domains.

A value of $k_{\text{ACL}}$ can be determined using the Darcy-Weißsbach equation, assuming a pore diameter of 52~nm \cite{alia_electrolyzer_2019} and a correction for the friction factor of 0.005 \cite{liu_experimental_2012}. An ACL thickness of 1~$\mu m$ is assumed for the ACL with low catalyst loading \cite{alia_electrolyzer_2019}. $P^{\text{atm}}$ is the atmospheric pressure which can be considered as the pressure level inside the pores of the PTL domain.
 
Equation \eqref{eq:balance_H} considers output and input flows. The outflow is defined as $F^{\text{PTL}}_{\ce{H2}}$ in equation \eqref{eq:ptl_flow_o2_h2}. The input flow $F^{\text{Mem}}_{\ce{H2}}$ which comes from the membrane domain will be explained in the following Subsection \ref{sec:1.2}.

Equation \eqref{eq:dissolution} describes the temporal evolution of the amount of Ir in the form of nanoparticles ($N_{Ir}^{np}$ measured in moles) within the ACL due to the two previously described dissolution mechanisms. The ECSA, defined by $A^{\text{act}}_{\ce{Ir}}$, is calculated based on
\begin{equation}
A^{\text{act}}_{\ce{Ir}} = {\eta}^{\text{act}}_{\text{np}} 4 \pi r_{\text{np}}^2 \quad \text{and} \quad \quad r_{\text{np}} = \sqrt[3]{\frac{3 N^{np}_{\ce{Ir}} M_{\ce{Ir}}}{4 \pi \rho_{\text{Ir}} {\eta}^{\text{act}}_{\text{np}}}}.
\end{equation}
The variable \({\eta}^{\text{act}}_{\text{np}}\) represents the number of active nanoparticles, while the radius \(r_{\text{np}}\) denotes their size. Additionally, parameters such as the molar mass of Iridium, denoted as $M_{\ce{Ir}}$ (in [g/mol]), and its density $\rho_{\text{Ir}}$ (in [g/m3]) are present in the formulation. 

In the kinetics a linearity between Ir catalyst dissolution and OER at higher current densities is included. This linearity results in the same factor in front of the potential within the exponent of the rate expression.

\subsection{Transport of hydrogen through the membrane}\label{subsec:transport_h2}\label{sec:1.2}
\begin{figure}[h]
    \begin{center}
	  \includegraphics[width=1.0\textwidth]{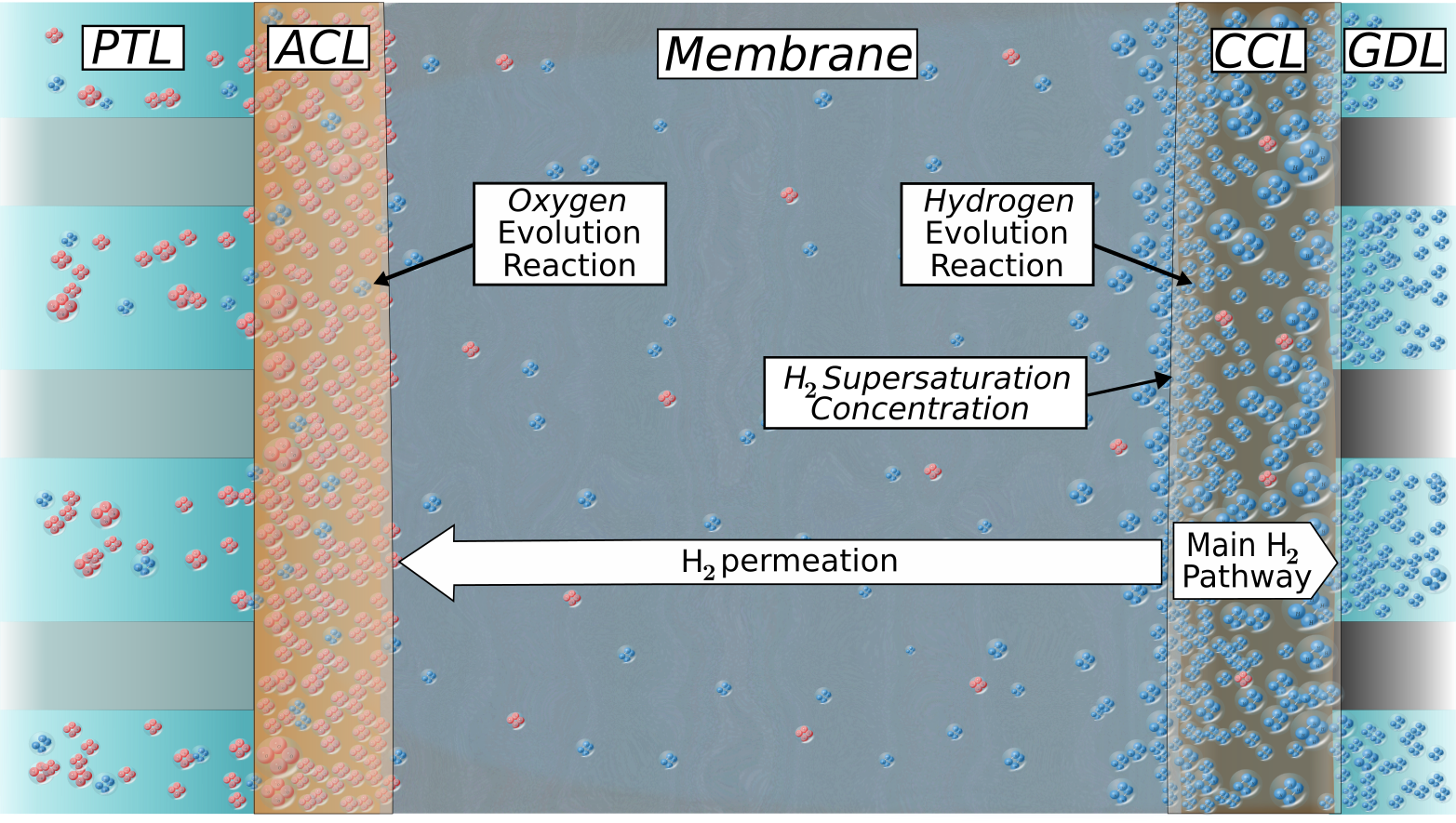}  
    \end{center}
    \caption{The sketch shows a section of the PEMWE consisting of the membrane at the center and the ACL and the CCL at the sides. In the cathode catalyst layer most of the evolved hydrogen is transported to the cathode outlet. However, a small fraction of hydrogen permeates into the other direction towards the anode side. In the ACL oxygen is evolved and a small volumetric fraction is constituted by hydrogen.}
    \label{fig:H2_transport}
\end{figure}
Due to the electrochemical polarization of the catalyst layer, hydrogen is evolved in the CCL. In previous studies, it has been found that this creates a supersaturation of hydrogen in the CCL, resulting not only in the transport of hydrogen in the cathode channel direction but also the permeation of hydrogen to the anode side \cite{trinke_current_2017}. This is schematically shown in Figure \ref{fig:H2_transport}. In this work, the transport is modelled via a Partial Differential Equation which considers the time-dependent diffusion of hydrogen over the membrane domain as
\begin{equation}
\partial_t C^{\text{MEM}}_{\ce{H2}}(x, t) - D^{\text{eff}}_{\ce{H2}}\partial_{xx} C^{\text{MEM}}_{\ce{H2}}(x, t) = 0.\label{eq:transport}
\end{equation}
The PDE is discretized and solved using the Finite Element Method in a one-dimensional domain, considering $D_{\ce{H2}}^{\text{eff}}$ as the effective hydrogen diffusion coefficient of the membrane.

The boundary condition at the CCL side was implemented as a Dirichlet boundary condition. For simplicity, this can be done by assuming quasi-steady-state conditions at the CCL and therefore calculating the supersaturation concentration based on the following equation (\cite{trinke_current_2017}):
\begin{equation}
C_{\text{\ce{H2}}}^{\text{CCL}}(t) = \frac{\frac{i}{2 F} + k_l C_{\text{\ce{H2}}}^{\text{Henry}}(t)}{k_l + \frac{D^{\text{eff}}_{\ce{H2}}}{\delta^{\text{MEM}}}}.\label{eq:boundary_conditions_CCL}
\end{equation}
On the rhs, the first term of the numerator represents the $\ce{H2}$ evolution which depends on the current density $i$. For the mass transfer coefficient $k_l$, which describes the transport from the reaction local to the cathode GDL domain, a value of $10^{-4} m/s$ is used, taken from the same source. $C_{\text{\ce{H2}}}^{\text{Henry}}$ is the theoretical saturation concentration of dissolved hydrogen assumed from Henry's law without mass transfer limitations. The membrane $\delta^{\text{MEM}}$ is considered one-dimensionally in space and discretized from $0$ to the membrane thickness of $L=183~\mu m$.

On the other side, the boundary conditions for the membrane and the ACL are defined via Neumann boundary conditions
\begin{equation}
D^{\text{eff}}_{\ce{H2}}\partial_x C^{\text{MEM}}_{\ce{H2}}(0, t) = k^{\text{MEM}}_{\ce{H2}} R T \left( C^{\text{MEM}}_{\ce{H2}}(0, t) - C^{\text{ACL}}_{\ce{H2}}(t)\right).
\end{equation}
For the mass transfer coefficient $k^{\text{MEM}}_{\ce{H2}}$ a very high value of $10^{-10}$ is chosen and therefore, effectively, the resistance directly at the interface of membrane and ACL domain is considered negligible. Note that the entire equation has been divided by the diffusion coefficient, which is consequently combined into the model parameter $k^{\text{MEM}}_{\ce{H2}}$ (unit $1/m$).

The factor $RT$ comes from the original definition of the transport rate coefficient, which is based on the pressure difference as the driving force. The ideal gas law, as given by equation \eqref{eq:idealGasLaw}, can be used to express the same term via the concentration difference. $C^{\text{MEM}}_{\ce{H2}}(0, t)$ is the estimated concentration of hydrogen at the first element of the discretized membrane over time and the concentration inside the ACL is given by the variable $C^{\text{ACL}}_{\ce{H2}}$. The flux between the membrane and the ACL  $F^{\text{Mem}}_{H_2}$ in equation \eqref{eq:balance_H} is calculated by:
\begin{equation}
F^{\text{Mem}}_{H_2}(t) = k^{\text{MEM}}_{\ce{H2}} R T ( C^{\text{MEM}}_{\ce{H2}}(0, t) - C^{\text{ACL}}_{\ce{H2}}(t) ).
\end{equation}
This relationship is employed to couple both problems: the ACL modelled by an ODE system and membrane domain by a PDE. The coupling is achieved through a Fixed Point Iteration Method as explained in Section \ref{sec:3}.

\section{Temporal Multiscale Method for solving the degradation model}
\label{sec:2}
Upon closer examination of the complete problem described by equations (\ref{eq:kineticks} - \ref{eq:dissolution}, \ref{eq:transport}) a difference in time scales between equation \eqref{eq:dissolution} and the others becomes noticeable. The process of dissolution, in fact, takes hours to undergo a measurable change in its state. In contrast, the remaining equations describe very fast processes as kinetics (represented by $\Theta_1$ and $\Theta_2$) and the changes in the concentrations of $\ce{O2}$ and $\ce{H2}$ in the ACL. As a matter of effect, those occur in time lapses less than a second.

The assumed magnitudes for the mentioned equations are responsible for the time scale distinctions. The first term of \eqref{eq:kineticks} is approximately $10^{-11}$ taking into account that $k_r \approx 10^{-14}$, $0 \leqslant \Theta_1 \leqslant 1$ and the maximum value of $C^{\text{ACL}}_{H_2}$ is approximately $0.5$ bar. The second term has a magnitude of $10^{-12}$ due to $k_{\text{diss}_1} \approx 10^{-25}$, resulting in a total magnitude of $10^{-12}$ for the entire equation. In the case of the equation \eqref{eq:dissolution}, the constant $k_{\text{diss}_2}$ is assumed to be around $10^{-37}$, resulting in a total magnitude of approximately $10^{-5}$.

Therefore, there exists a proportionality between the variables which can be approximated by $\varepsilon \approx 10^{-7}$. This reveals an underlying temporal multiscale structure of the complete problem, as proposed by Weinan in \cite{Weinan2011}, with variables categorized into Fast-Scale and Slow-Scale from equations (\ref{eq:kineticks} - \ref{eq:balance_H}, \ref{eq:transport}) and \eqref{eq:dissolution} respectively. Furthermore, the Slow-Scale is $O(1)$ with respect to the Fast-Scale through $\varepsilon$.

Moreover, the Fast-Scale processes also display a local quasi-periodic behavior due to the input of the potential. If the slow variable would not further evolve, a periodic-in-time forcing would result in a periodic response. Thus, these profiles simulate the oscillatory behavior resulting from real-world degradation of PEMWEs in operation. Then, during a $P$ period of time there are not significant changes of the Slow-Scale variables.

A direct simulation of such a temporal multiscale problem is not efficient. Resolving the complete time interval with time steps of less than a second (e.g. $\Delta k = 10^{-2}s$) would result in excessive computational effort although the slow process, which is the one of most interest, is hardly changing. We therefore utilize a temporal multiscale algorithm that is based on the heterogeneous multiscale method \cite{Weinan2003,Abdulle2012} and that aims at decoupling slow and fast scale. 

The first step of decoupling is to derive an effective equation by averaging the slow process \eqref{eq:dissolution} over one period $P$. Second, instead of resolving the fast process dynamically in time, it is localised by fixing the slow variable $N^{np}_{Ir}$ and by replacing the Fast-Scale feedback by the periodic-in-time limit cycle of (\ref{eq:kineticks} - \ref{eq:balance_H}, \ref{eq:transport}) which can be numerically approximated without simulating the complete history of the fast scale variables. By this decoupling, a temporal multiscale method is obtained that allows to step forward the coupled temporal multiscale problem with very large time-steps. We refer to \cite{frei_efficient_2020,Lautsch2021} for details.  

\begin{figure}[h]	
    \begin{center}
		\includegraphics[width=1\textwidth]{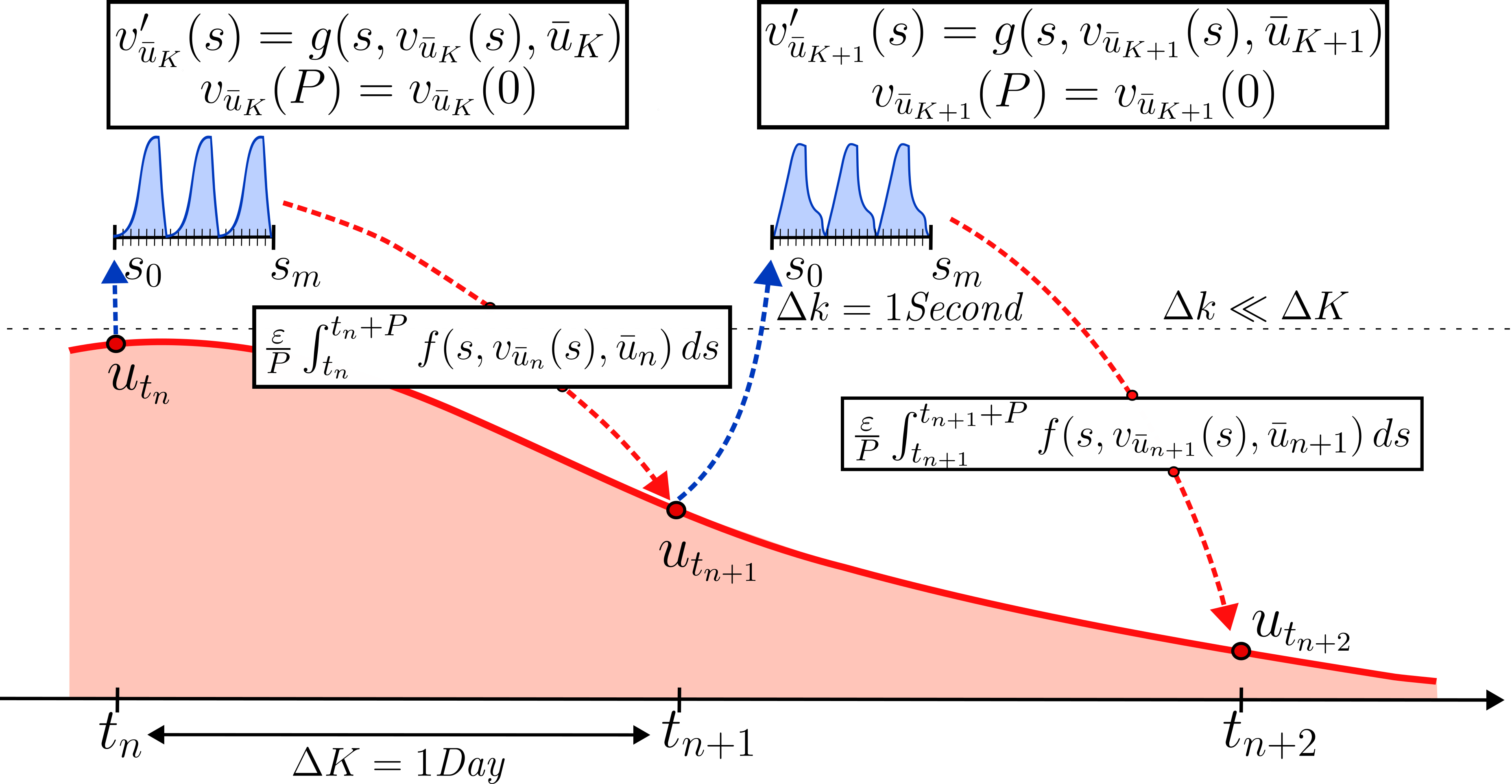}
    \end{center}
    \caption{Schema of the temporal multiscale method.}
    \label{fig:multirate_method}
\end{figure}

Taking advantage of this particular structure of the problem and the theory of Averaging of the multiscale problems it is possible to approximate the complete problem with its corresponding temporal multiscale problem as in the next subsection \ref{sec:2.1} will be explained.

\subsection{Temporal Multiscale Method}\label{sec:2.1}
Based on the stated, the Slow-Scale variables of the approximated temporal multiscale problem are defined over $t \in [0, T] \subset \mathbb{R}$ having $\Delta K$ as a large time step size. In the case of the Fast-Scale variables, if $P$ is established as the period and $0 \leqslant m \leqslant (T - \Delta K) / \Delta K$, then these variables will be defined over $s \in [t_K + \Delta K + n P, t_K + \Delta K + (n+1) P]$ for $0 \leqslant n \leqslant (\Delta K - P) / P$  and $t_K = m \Delta K$ with a small time step size $\Delta k \ll \Delta K$. The new problem is written as follows:

Find $\bar{N}^{np}_{Ir}\left( t\right)$ that approximates $N^{np}_{Ir}\left( t\right)$ and is the solution of 
\begin{equation}
\dfrac{d \bar{N}^{np}_{Ir}}{dt} = -\frac{\varepsilon}{P} \int_{t_{K} + n P}^{t_{K} + (n+1) P} \frac{A^{\text{act}}_{\ce{Ir}}}{A_{geo}} \left[\hat{\Theta}_1 \left( s\right) k_{\text{diss}_1} + \left( 1 - \hat{\Theta}_1\left( s\right) \right) k_{\text{diss}_2} \right] e^{f E^{\text{ACL}} \left( s\right)} \, ds\label{eq:ms_dissolution}
\end{equation}
so that for all $K$ there exists $n \in \mathbb{N}_0$ and the interval $[t_{K} + n P, t_{K} + (n+1) P]\subset [t_{K}, t_{K+1}]$ such that the equations of the Fast-Scale
\begin{align}
\dfrac{d \hat{\Theta_1}\left( s\right)}{ds} &= 2 k_r (1- \hat{\Theta}_1\left( s\right))^2 \left( \hat{C}^{\text{ACL}}_{H_2}\left( s\right) R T \right) - \frac{\hat{\Theta}_1\left( s\right)}{\gamma} k_{\text{diss}_1} e^{f E^{\text{ACL}}\left( s\right)}\label{eq:ms_kineticks}\\
\dfrac{d \hat{C}^{\text{ACL}}_{O_2}\left( s\right)}{ds} &= \frac{A^{\text{act}}_{Ir}\left(\bar{N}^{np}_{Ir} \left( t \right)\right)}{V_{\text{ACL}}}\left[ \frac{i}{z_{O_2} F} - F^{PTL}_{O_2}\left( s\right) \right]\label{eq:ms_balance_O}\\
\dfrac{d \hat{C}^{\text{ACL}}_{H_2}\left( s\right)}{ds} &= \frac{A_{\text{geo}}}{V_{\text{ACL}}}\left[ F^{\text{Mem}}_{H_2}\left( s\right) - F^{\text{PTL}}_{H_2}\left( s\right)\right]\label{eq:ms_balance_H}\\
\partial_s \hat{C}^{\text{MEM}}_{\ce{H2}}&(x, s) = D^{\text{eff}}_{\ce{H2}}\partial_{xx} \hat{C}^{\text{MEM}}_{\ce{H2}}(x, s),\label{eq:ms_transport}
\end{align}
with the same initial conditions as the complete problem, are satisfied in each chosen interval and their solution is locally quasi-periodic. 

The choice of $\Delta K$ is conditioned by the prior choice of $\Delta k$ and the previously calculated scaling constant $\varepsilon$, such that $\Delta K \approx \Delta k / \varepsilon$. In this case, $\Delta k$ is set to $10^{-2}$ seconds to ensure stability over the Fast-Scale resolution. Since 	the experimental data used for validation was sampled in days, it was convenient to choose $\Delta K = 1$ day, which is approximately $10^5$ seconds and satisfies the mentioned temporal scale relation.

In order to simplify the notation some new symbols have been included. The $(\bar{.})$ accent will be understood as the average of the function, in this case the Slow-Scale variable. On the other hand, $(\hat{.})$ will represent the dependency 
\begin{equation}
\hat{v} = (\hat{\Theta}_1, \hat{C}^{\text{ACL}}_{O_2}, \hat{C}^{\text{ACL}}_{H_2}, \hat{C}^{\text{MEM}}_{\ce{H2}})
\end{equation} 
in which function $v$ could be any of the Fast-Scale variables (in case of variable $\hat{C}^{\text{MEM}}_{\ce{H2}}$, $x$ is also included).

\section{Numerical implementation}\label{sec:3}
The numerical solution of the proposed problem in Section \ref{sec:2} is achieved by the implementation of the general temporal multiscale method shown in Algorithm \ref{alg:MSM}.
\begin{algorithm}
\caption{General temporal multiscale method for PEMWE dissolution problem}\label{alg:MSM}
\begin{algorithmic}[1]
\Require $ u_K = \bar{N}^{np}_{Ir}\left( t_0\right), \hat{v}_{u_K} = v_0$
	\For{ $K = 0, \dots, N-1$}  
		\State \textit{ Find the solution of the Fast-Scale equations $\hat{v}_{u_K}$ with fixed $u_K$ such that for $n>0$:}
		\Statex $||\hat{v}_{u_K}(t_K + n P) - \hat{v}_{u_K}(t_K + (n+1) P)||< tolp$
		\State \textit{Find $u_{K+1}$ integrating the average of the Slow-Scale equation in the interval of the Fast-Scale solutions $\hat{v}_{u_K}$ such that:}
		\Statex $u_{K+1} = u_K - \frac{\varepsilon \Delta K}{P} \int_{t_{K} + n P}^{t_{K} + (n+1) P} \left[\hat{\Theta}_1 \left( s\right) k_{\text{diss}_1} + \left( 1 - \hat{\Theta}_1\left( s\right) \right) k_{\text{diss}_2} \right] e^{f E^{\text{ACL}} \left( s\right)} \, ds$
	\EndFor\\
	\textbf{return} $u_{K+1}$
\end{algorithmic}
\end{algorithm}
The main idea is to solve the Fast-Scale problem using a \textit{semi-implicit} method with a fixed $u_K$ solution over the interval $[t_K, t_K + \Delta K]$, until it reaches a locally quasi-periodic solution $\hat{v}_{u_K}$ with an error tolerance $||tolp||_{\infty}$. Then, the periodicity convergence criteria is defined as the infinite norm of the component-wise absolute errors for each Fast-Scale variable of the ACL. After that, $u_{K+1}$ is obtained by integrating the averaged Slow-Scale equation over the time interval which $\hat{v}_{u_K}$ is defined. 

In order to solve the Fast-Scale problem the ODE System (\ref{eq:ms_kineticks} - \ref{eq:ms_balance_H}) coupled to the Transport Problem \eqref{eq:ms_transport} must be solved. The solution of the Transport Problem was achieved by the application of a one dimensional time implicit Finite Element Method. Approximation of integrals of the mass and stiffness matrices were carried out by applying the trapezoidal method. For time integration the implicit backward Euler method was implied. 

On the other hand, due to the stiff features of the problem, the ODE System mentioned before was solved applying an second-order backward differentiation formula scheme implemented in the `ode' function of the `integrate' module of the scientific computing library `scipy' \cite{2020SciPy-NMeth}, available in the Python programming language. Back in Section \ref{sec:1}, it was explained how hydrogen coming from the membrane could drastically permeate the ACL, and this is precisely what causes a jump in the solutions, making it more difficult to solve with explicit methods, which will need time step adaptivity and therefore longer computation time.

Coupling both problems in each time step $\Delta k$ was accomplished by first solving the ODE system (\ref{eq:ms_kineticks} - \ref{eq:ms_balance_H}) implicitly with step size $\Delta k$, keeping $u_K$ and $\hat{C}^{\text{MEM}}_{\ce{H2}}(t_K)$ fixed. Afterwards, the Transport Problem \eqref{eq:ms_transport} is solved considering the solutions of the system (\ref{eq:ms_kineticks} - \ref{eq:ms_balance_H}) corresponding to $t_K$ and $t_{K+1}$. Once a locally quasi-periodic solution $\hat{v}_{u_k}$ is found, then the Slow-Scale equation is integrated. The integration method used for this purpose is the Adams-Bashforth second order method.

\section{Results and discussion}\label{sec:4}
\begin{figure}[h]
	\begin{center}
		\includegraphics[width=1\textwidth]{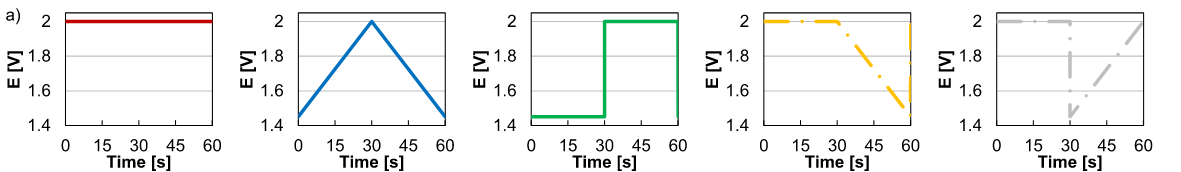}
		\caption{Experimental data of five test profiles acquired from \cite{alia_electrolyzer_2019}. They are named from left to right: hold, triangle-wave, square-wave, sawtooth-up and sawtooth-down profiles.}
		\label{fig:Alia_profiles}
	\end{center}
\end{figure}

The dynamic operation profiles, as well as experimental data, are taken from Alia et al. \cite{alia_electrolyzer_2019}. The cyclic voltammetric charge, which was systematically measured in this work, is used as a metric for the electrochemical surface area of the PEMWE device. In this experimental investigation, five dynamic operation profiles were applied in which the cell potential was controlled (Figure \ref{fig:Alia_profiles}). These long-term experiments were conducted for 525 hours of operation.

\subsection{Temporal Multiscale simulation}
A comparison between the fully-resolved and the temporal multiscale problems was made to verify the accuracy of the approximation that could be achieved with the multiscale approach. The simulations where ran simulating 21 days of operation of the PEMWE with initial ECSA $= 3.25 \cdot 10^{-2} m^2$, period of $60$ seconds and, as was shown in Section \ref{sec:2.1}, $\Delta k$ and $\Delta K$ are chosen as $10^{-2}$ seconds and $24$ hours respectively.  
\begin{table}[ht]
	\caption{Computational time and mean squared error (MSE) of the results for the fully resolved problem and the temporal multiscale method (FRP and TMS) across the five experimental profiles \cite{alia_electrolyzer_2019} for 21 simulated days.}
    \centering
    \begin{tabular}{ccccc}
        \toprule
        \multirow{2}{*}{Experiment} & \multicolumn{2}{c}{Computational time} & \multirow{2}{*}{$\text{MSE}\big(N_{\text{Ir}}^{\text{np}}, \bar{N}_{\text{Ir}}^{\text{np}} \big)$} & \\
        \cmidrule(r){2-3}
          profiles & FRP(days) & TMS(min) & \\
        \midrule
        Hold & 8.6 & 19.09 & $7.87 \cdot 10^{-06}$ \\
        Square-wave & 6.8 & 13.29 & $2.46 \cdot 10^{-04}$ \\
        Triangle-wave & 6.4  & 13.46  & $5.39 \cdot 10^{-06}$ \\
        Sawtooth Down & 8.0 & 16.84 & $3.00 \cdot 10^{-05}$ \\
        Sawtooth Up & 8.2 & 21.06 & $5.52 \cdot 10^{-04}$\\
        \bottomrule
    \end{tabular}    
    \label{tab:experiment_data}
\end{table}

In Table \ref{tab:experiment_data}, the results of this comparison, considering the already fitted parameters, are tabulated. From it, we see that the average ratio of computational time between these two methods is approximately $2.82$ orders of magnitude, which means that the temporal multiscale method is around thousand times faster than the fully-resolved one on average. The best computational time for the fully-resolved problem is about $6.4$ days, which, compared with the $21$ minutes that the longest temporal multiscale simulation took, represents a significant computational reduction. Thus, approximately $1.25$ simulation days are computed per $1$ minute of computation for the temporal multiscale method. 

The maximum mean squared error between these two sets of simulation results is approximately  $5.52 \cdot 10^{-4}$, which is very small considering the enormous computational reduction made. This shows that the temporal multiscale method is accurate enough to be used instead of the fully-resolved one for approximating the long-term effects of the dissolution of Iridium when a PEMWE is operated. This property is ideal for estimating parameters efficiently as will be shown in Section \ref{sec:4.2}.

\subsection{Parameter fitting}\label{sec:4.2}
Regarding fitting the parameters of the model, a single set of parameters $k_r$ and $k_{\text{diss}_2}$, was estimated using data obtained from hold and square-wave profiles. The data from the other profiles was used as a test set. The objective function of the resulting optimization problem is expressed as
\begin{equation}\label{eq:fitting_eq}
\sigma^* = \arg\min_{\sigma \in \Omega \subset \mathbb{R}^2} \frac{1}{T}\int_0^T(S_{\text{HP}}(t) - u_{\text{HP}}(t; \sigma))^2 + (S_{\text{SWP}}(t) - u_{\text{SWP}}(t; \sigma))^2 \, dt,
\end{equation} 
where $\Omega$ is the space of feasible solutions of the form $\sigma = (k_r, k_{\text{diss}_2})$, $S_{\text{HP}}$ and $S_{\text{SWP}}$ the sampled data from hold and square-wave profile and $u$ the Slow-Scale variable for each case.

For the fitting procedure, a variant of the Simulating Annealing Optimization Method called Simple Serializable Set from \cite{GREENING1990293} was applied. In this way, a computer server with 128 cores and 260 GB of memory ram was used to widely spread the search around partial solutions, running 128 simultaneous simulations for each set of parameters in the search neighborhood for each iteration of the method. At the end of the fitting procedure the values $k_r^* = 2.437792 \cdot 10^{-14}$ and $k_{\text{diss}_2}^* = 1.947001 \cdot10^{-37}$ have been identified. The mean squared errors of Table \ref{tab:experiment_data} show the achieved accuracy.

\subsection{Evaluation of the Model Performance with Experimental Data} 
\begin{figure}[h]
    \centering
	\includegraphics[width=1.03\textwidth]{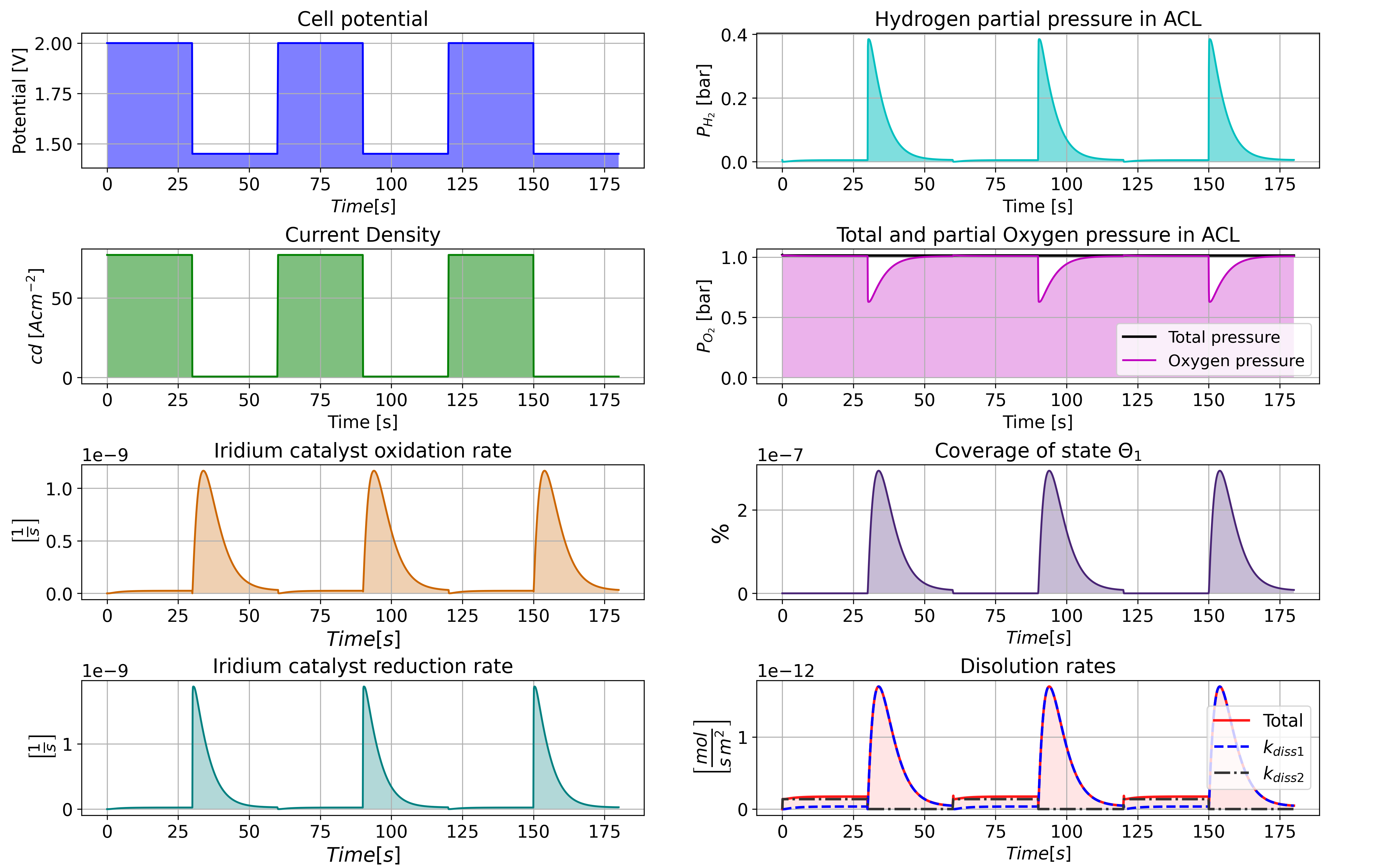}
    \caption{Fast-Scale variables for the square-wave profile simulation for the first 3 periods ($0 \leq t \leq 180$ seconds). The solution of the transport of hydrogen through the membrane is shown in the }\label{fig:Fast_Scale_Output}
\end{figure}
\begin{figure}[h]
    \centering
	\includegraphics[width=1\textwidth]{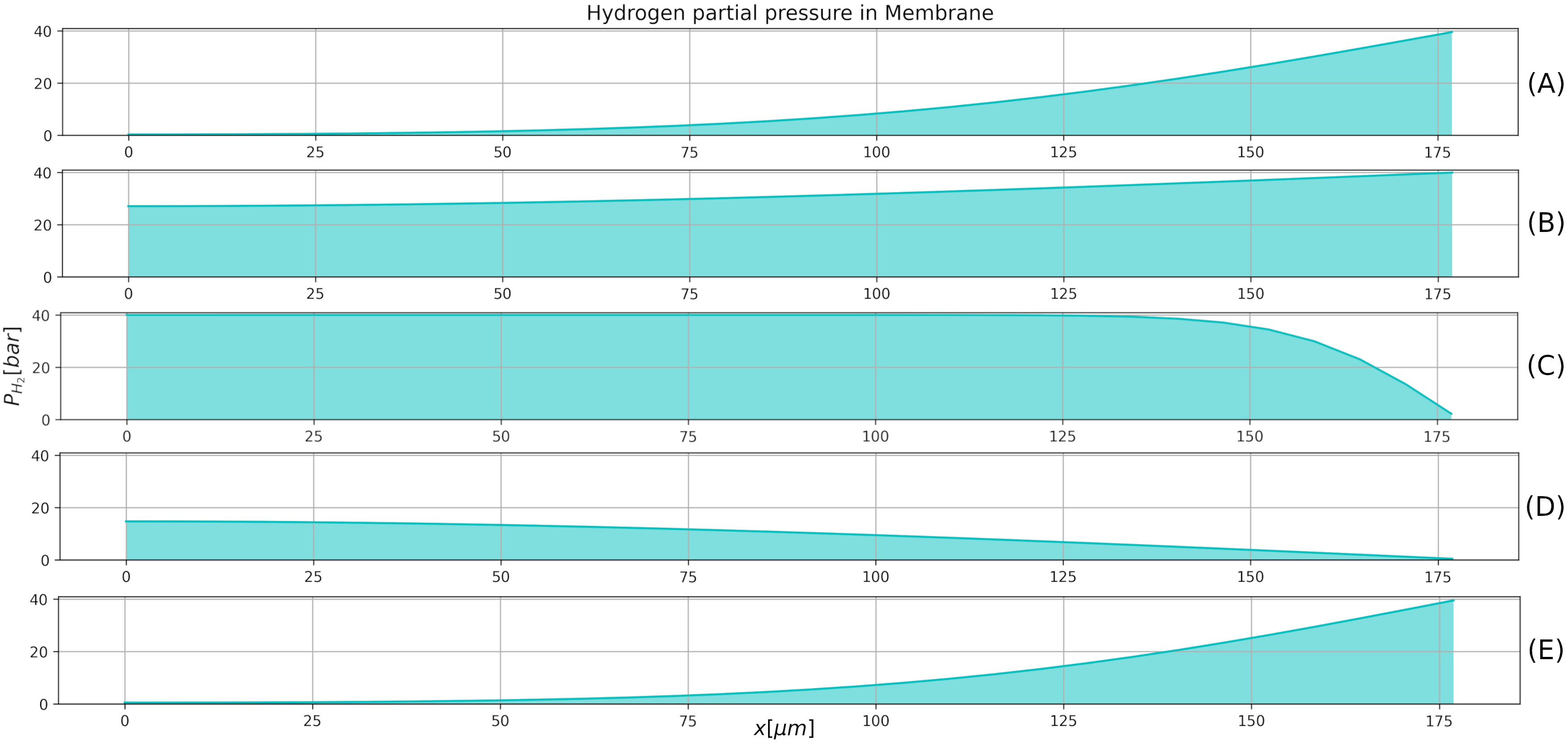}
    \caption{The Fast-scale solutions of the hydrogen transport problem in the membrane for the period $t\in[120, 180]$ sec. The graph (a) is the state of the concentration of hydrogen at $t=120$ sec, (b) at $t=125$ sec, (c) at $t=150$ sec, (d) at $t=155$ sec and (e) at $t=180$ sec respectively.}\label{fig:Transport_Problem_PDE}
\end{figure}

Multiple output variables were obtained for each Fast-Scale variable. The Fast-Scale variables of the hold profile exhibit steady-state behavior almost throughout the entire duration. To provide a more insightful comparison, the results of the square-wave profile are presented in Figure \ref{fig:Fast_Scale_Output} and Figure \ref{fig:Transport_Problem_PDE}, as this profile, along with the Hold profile, was chosen for fitting the model parameters.

In Figure \ref{fig:Fast_Scale_Output}, abrupt changes in all variables are observable at every period, specifically at 30 seconds after the start of each. These changes result from sudden drops in the applied cell potential, which is also the source of the \textit{stiffness} discussed in Sections \ref{sec:1} and \ref{sec:3}.

With respect to the transport of hydrogen in the membrane, Figure \ref{fig:Transport_Problem_PDE} illustrates various stages of hydrogen concentration profile during the last period shown in Figure \ref{fig:Fast_Scale_Output}. Graphs (a) and (b) indicate the build up of the pressure during the rise of potential up to $1.45$V.  An initial stage of potential drop is shown in (c), where the hydrogen concentration decreases in the CCL due to the drop in current density. In (d), the amount of hydrogen stored in the membrane has been almost completely transported out. Finally, in (e) a high current density was restarted, and the entire process repeats itself.

\begin{figure*}[h]
    \centering
    \begin{subfigure}[t]{0.48\textwidth}
        \centering
        \includegraphics[width=1.05\textwidth]{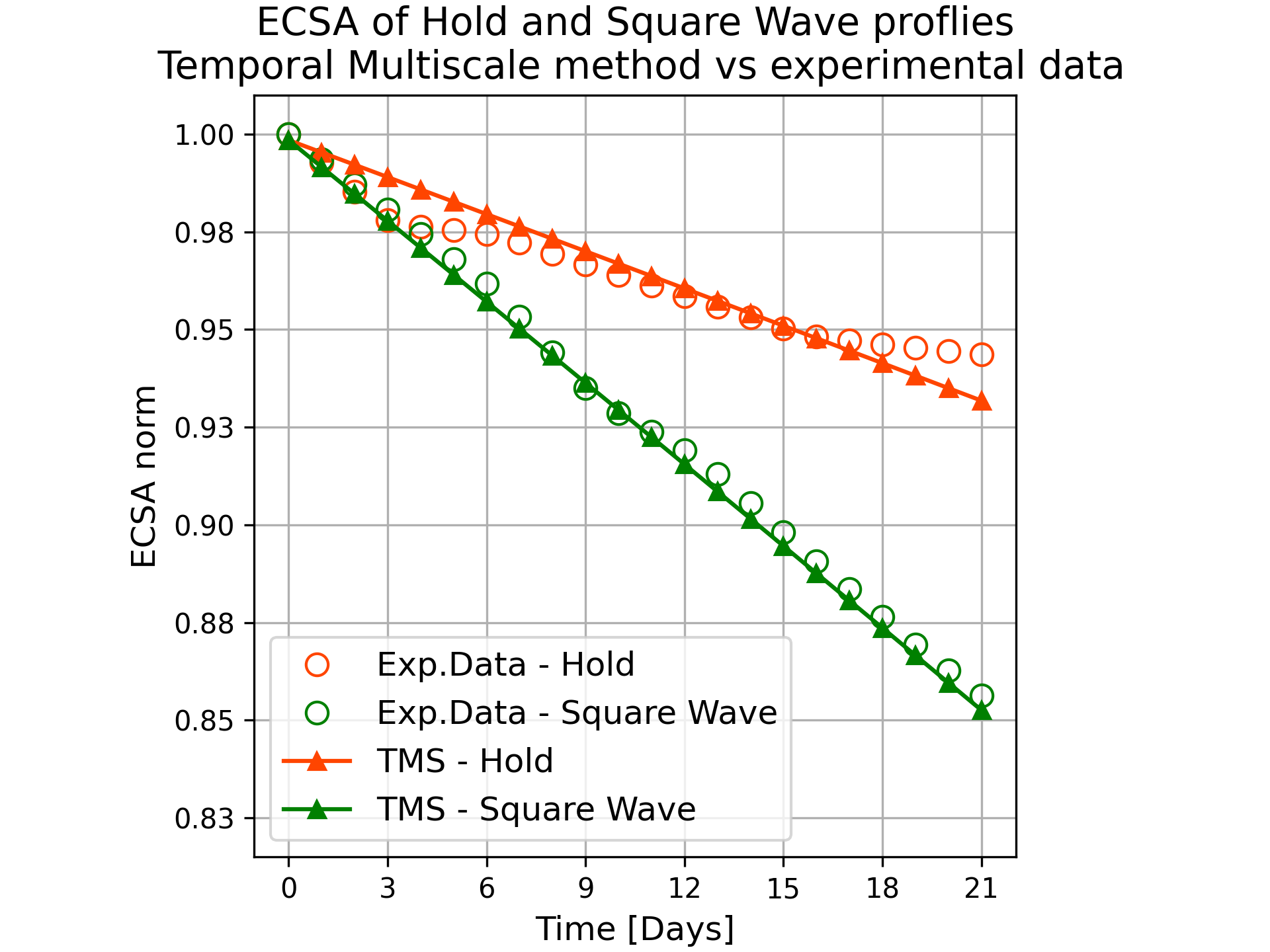}
        \caption{Parameter estimation step: Normalized ECSA data from Alia \cite{alia_electrolyzer_2019} vs the simulated results from the model with the application of the temporal multiscale method (TMS) for hold and square-wave profiles.}\label{fig:SlowScaleResults}
    \end{subfigure}%
    ~
    \hspace{0.02\textwidth}
    ~
    \begin{subfigure}[t]{0.48\textwidth}
        \centering
        \includegraphics[width=1.05\textwidth]{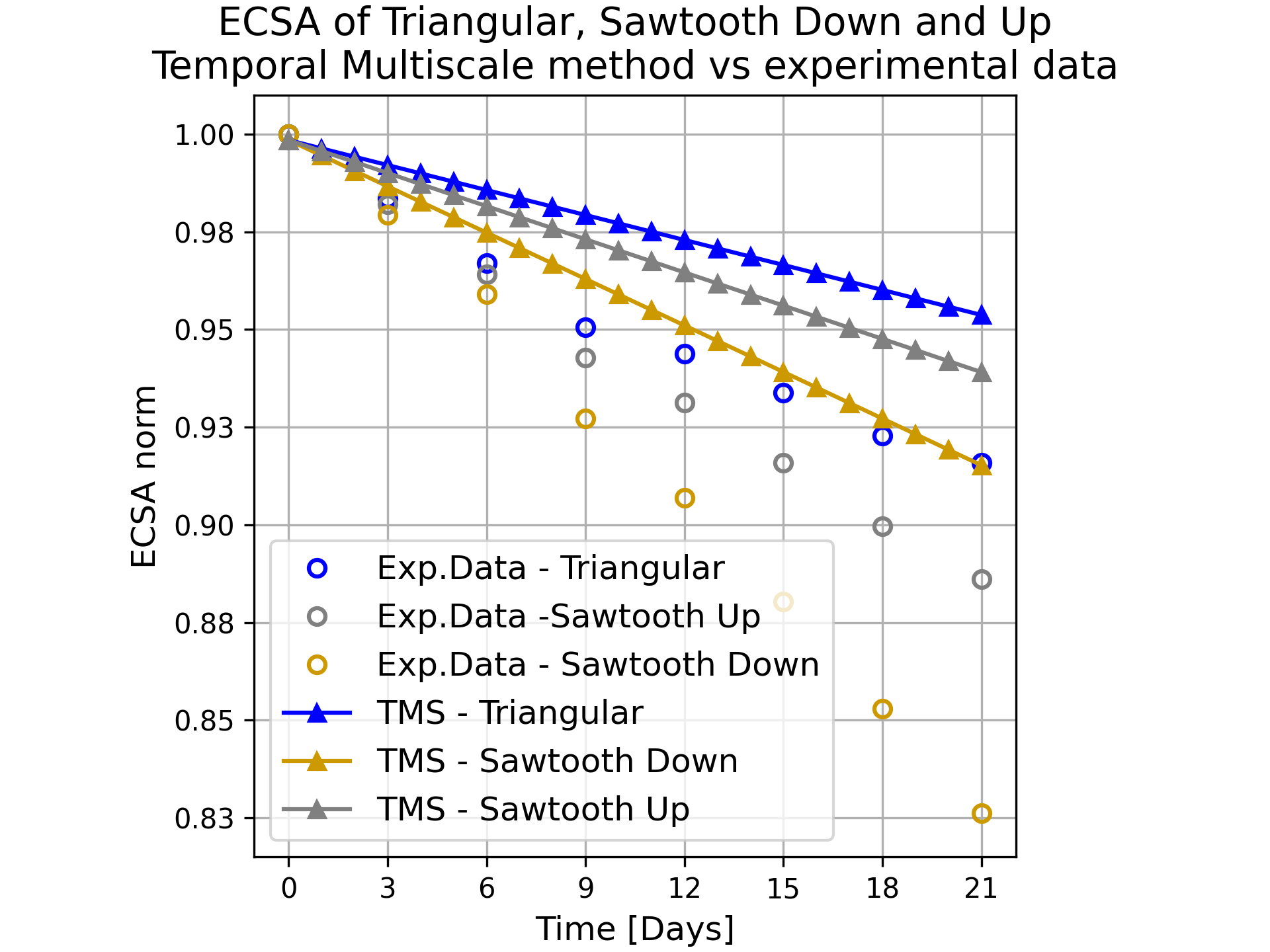}
        \caption{Model prediction step: Normalized ECSA data from Alia \cite{alia_electrolyzer_2019} vs the simulated results from the model with the application of the temporal multiscale method (TMS) for triangular, sawtooth-up and down profiles.}
    \end{subfigure}    
    \caption{Slow-Scale variables results compared with the experimental data provided by \cite{alia_electrolyzer_2019} after 21 days of simulation.}\label{fig:SlowScaleResults} 
\end{figure*}
The Slow-Scale variable is presented in Figure \ref{fig:SlowScaleResults}. It is illustrated how the ECSA decays over time for experimental data and the model simulation. It can be observed that, through the described fitting procedure, a single set of parameters has been obtained. The use of the latter results in a good agreement between the simulation and experimental data for both the hold and square-wave operation profiles. Only a slight bending up observed in the experimental data of the hold profile does not seem to be well captured by the implemented model. This may indicate that a mechanism  leading to the stabilization of the catalyst under steady-state conditions of the hold profile has been neglected by the model. However, apart from this aspect, based on the minimized low root mean squared error, a successful model fitting and therefore parameter estimation was achieved. 

This set of parameters is then employed to predict the outcome of the other three PEMWE operation profiles that are presented in Figure \ref{fig:Alia_profiles} (triangle, sawtooth-up and sawtooth-down). The model simulations and the corresponding data are presented in Figure \ref{fig:SlowScaleResults} b). The model predictions show a significant deviation from the experimental data, indicating The model predictions show a significant deviation from the experimental data, indicating the need for further development and improvement of the mechanistic degradation model. In particular, the model appears to underestimate the degradation effect of the slower, dynamic ramps of the triangular and the two sawtooth operation profiles.

Finding the `right' electrochemical model which fits all the dynamic operation profiles shown is out of the scope of this work. However, it is shown how the application of the temporal multiscale approach can lead to a much faster simulation of long-term degradation processes in PEMWE. 

A systematic model development can be performed by iterating the following steps: 1) proposing/modifying the mathematical degradation model, 2) applying the parameter fitting procedure, 3) performing predictive simulations and 4) comparing the predictive simulations with the experimental data. The temporal multiscale approach can play an important role in this iterative process by facilitating the simulation and especially the parameter estimation step. Once a reliable degradation model is found, the temporal multiscale approach can also be used to optimize the dynamic operation of PEMWE, taking into account the long-term degradation effects of different dynamic operation strategies.

\section*{Conclusion and future perspectives}
In this work, a temporal multiscale method was applied to accelerate simulation of degradation in PEMWE under dynamic operating conditions. With an exemplary focus on the phenomenon of catalyst dissolution, a mechanistic degradation model was hypothesized and implemented. Due to the specific underlying temporal multiscale structure of the fully-resolved problem, its computational complexity was reduced by solving a temporal multiscale equivalent problem. 

In this regard, two different temporal scales were identified: A fast process with quasi-local periodicity dictated by the transport and reaction processes occurring under intermittent operation conditions of PEMWE, and a slow process constituted by the gradual degradation of the catalyst layer.
 
Assuming that the Slow-Scale variables are minimally affected by the Fast-Scale variables in a time step significantly larger than the one used for the time discretization of the Fast-Scale variables, we can consider the Slow-Scale variables as fixed during the computation of the Fast-Scale solution. Additionally, the presented temporal multiscale method leverages the quasi-local periodicity of the Fast-Scale variables not computing all the time steps but rather the necessary ones to achieve such periodicity. This enables the integration of the Slow-Scale variable solution over the mentioned large time step, utilizing the Fast-Scale variables solutions for just one period of time to compute the gradient.

The computational effort is substantially reduced compared to the fully-resolved problem solution. The temporal multiscale method, as mentioned before, doesn't required computing all the Fast-Scale solutions but only up to reaching a quasi-local periodicity for each Slow-Scale large time step. This makes it possible to reduce simulations of 21 days from 8 computational days for the fully-resolved problem to approximately 20 minutes. Furthermore, the mean squared error committed by the temporal multiscale method is approximately $10^{-5}$, making the impact of the assumptions negligible in terms of accuracy.

Numerical challenges, such as the stiffness of the Fast-Scale problem, were overcome by applying an implicit integration method using the Newton method. A semi-implicit integration method was also employed to ensure the stability of the coupling between the hydrogen transport and the ACL degradation problem.

The model parameters were fitted using a Simulating Annealing algorithm. Although a gradient descent method was not utilized as the optimization method for the parameter fitting  to avoid additional errors from applying a finite difference method for small values, the temporal multiscale approach proposed in this work also contributed to accelerating this procedure. Since the simulations took significantly less time than before, a wider search space of possible fitting solutions was explored. For this purpose, 128 simultaneous simulations were run for each iteration of the Simulating Annealing Method, resulting in a better approximation to the optimal parameters for fitting this problem. This could be replaced in the future by a fitting method that directly considers the particular temporal multiscale structure of the problem. 

The Slow-Scale variables were investigated to analyse the electrochemical problem.  The model parameters were estimated by fitting the model to the experimental data of two operation profiles (‘hold’ and ‘square-wave’) in PEMWE, resulting in a close agreement between the model simulation and the experimental data. However, the model prediction of the three other operation profiles ('triangular', 'sawtooth-up' and 'sawtooth-down') shows significant deviations from the experimental data, indicating the need for futher modification of the mathematical degradation model. 

In this paper, it is shown that the application of the temporal multiscale method can support the systematic development of a degradation model in PEMWE by drastically reducing the computational time for the model fitting and predictive simulation steps. In the future, the application of the temporal multiscale method can also aid to realize efficient numerical methods that can be used to optimize the dynamic operation of PEMWE.


\section*{Acknowledgements}
The work of DD was supported by the Center of Dynamic Systems (CDS), funded by the EU-programme ERDF (European Regional Development Fund).

In addition to this support, this work was authored in part by the National Renewable Energy Laboratory, operated by Alliance for Sustainable Energy, LLC, for the U.S. Department of Energy (DOE) under Contract No. DE-AC36-08GO28308. The views expressed in the article do not necessarily represent the views of the DOE or the U.S. Government. The U.S. Government retains and the publisher, by accepting the article for publication, acknowledges that the U.S. Government retains a nonexclusive, paid-up, irrevocable, worldwide license topublish or reproduce the published form of this work, or allow others to do so, for U.S. Government purposes.

\bibliographystyle{elsarticle-num-names} 
\bibliography{references}

\end{document}